\documentclass[a4paper]{article}
\usepackage{amsthm,amssymb,amsmath,enumerate,graphicx}
\advance\textheight by -1.6cm% to make page fit on my monitor at 2000

\def\proclaim #1.#2 #3\par{\bigbreak
   \noindent{\bf#1.}#2\enspace{\em#3}\par\bigbreak}

\newtheorem{definition}{Definition}

\newtheorem{theorem}[definition]{Theorem}
\newtheorem{corollary}[definition]{Corollary}
\newtheorem{lemma}[definition]{Lemma}

\newtheorem{problem}[definition]{Problem}

\newcommand{\COMMENT}[1]{}
\newcommand{\<}[1]{\vadjust{\vbox to 0pt{\vss\vskip-8pt\leftline{%
     \llap{\rm\hbox{\vbox{\pretolerance=-1
     \doublehyphendemerits=0\finalhyphendemerits=0
     \hsize16truemm\tolerance=10000\small
     \lineskip=0pt\lineskiplimit=0pt
     \rightskip=0pt plus16truemm\baselineskip8pt\noindent
     \hskip0pt        %(without this, the first word is never hyphenated!)
     #1\endgraf}\hskip7truemm}}}\vss}}}
\renewcommand{\<}[1]{}

\newcommand{\sm}{\smallsetminus}
\newcommand{\sub}{\subseteq}
\newcommand{\subne}{\subsetneq}
\newcommand{\supe}{\supseteq}
\newcommand{\es}{\emptyset}
\newcommand{\B}{\ensuremath\mathcal{B}}
\newcommand{\C}{\ensuremath\mathcal{C}}
\newcommand{\E}{\ensuremath\mathcal{E}}
\newcommand{\Emax}{\ensuremath\mathcal{E^{\rm max}}}
\newcommand{\I}{\ensuremath\mathcal{I}}
\newcommand{\Imax}{\ensuremath\mathcal{I^{\rm max}}}
\newcommand{\N}{\mathbb N}

\newcommand{\dcl}[1]{\ensuremath\lceil#1\rceil}% down-closure
\newcommand{\cl}{\ensuremath{{\rm cl}}}% closure operator

\newcommand{\mfc}{\ensuremath{M_{\rm FC}}}% finite cycles
\newcommand{\mc}{\ensuremath{M_{\rm C}}}% topological cycles
\newcommand{\mfb}{\ensuremath{M_{\rm FB}}}% finite bonds
\newcommand{\mb}{\ensuremath{M_{\rm B}}}% bonds
\newcommand{\atc}{\ensuremath{M_{\rm AC}}}% algebraic cycles

\newcommand{\emtext}[1]{\text{\em #1}}
\newcommand{\etop}[1]{\ensuremath{\|#1\|}}
\newcommand{\eends}{\mathcal E}
\newcommand{\kreis}[1]{\mathaccent"7017\relax #1}
\newcommand{\selfcite}[1]{$\!\!${\bf\cite{#1}}}

\newcommand{\xsided}{skew}%{small-sided}

\title{Infinite matroids in graphs}
\author{Henning Bruhn \and Reinhard Diestel}
\date{}

\begin{document}
\maketitle

\begin{abstract}\noindent
It has recently been shown that infinite matroids can be axiomatized in a way that is very similar to finite matroids and permits duality. This was previously thought impossible, since finitary infinite matroids must have non-finitary duals.

In this paper we illustrate the new theory by exhibiting its implications for the cycle and bond matroids of infinite graphs. We also describe their algebraic cycle matroids, those whose circuits are the finite cycles and double rays, and determine their duals. Finally, we give a sufficient condition for a matroid to be representable in a sense adapted to infinite matroids. Which graphic matroids are representable in this sense remains an open question.
\end{abstract}

\section{Introduction}

In the current literature on matroids, infinite matroids are usually ignored, or else defined like finite ones%
   \footnote{The augmentation axiom is required only for finite sets: given independent sets $I,I'$ with $|I| < |I'| < \infty$, there is an $x\in I'\sm I$ such that $I+x$ is again independent.} with the following additional axiom:

\begin{itemize}
\item[\rm (I4)] An infinite set is independent as soon as all its finite subsets are independent.
\end{itemize}

\noindent
We shall call such set systems {\em finitary matroids\/}.

The additional axiom (I4) reflects the notion of linear independence in vector spaces, and also the absence of (finite) circuits from a set of edges in a graph. More generally, it is a direct consequence of (I4) that circuits, defined as minimal dependent sets, are finite.

An important and regrettable consequence of the additional axiom (I4) is that it spoils duality, one of the key features of finite matroid theory. For example, the cocircuits of an infinite uniform matroid of rank~$k$ would be the sets missing exactly $k-1$ points; since these sets are infinite, however, they cannot be the circuits of another finitary matroid. Similarly, every bond of an infinite graph would be a circuit in any dual of its cycle matroid---a set of edges minimal with the property of containing an edge from every spanning tree---%
   \COMMENT{}%
   but these sets can be infinite and hence will not be the circuits of a finitary matroid.

This situation prompted Rado in 1966 to ask for the development of a theory of non-finitary infinite matroids with duality~\cite[Problem P531]{Rado66matroids}.%
   \COMMENT{}
In the late 1960s and 70s, a number of such theories were proposed; see~\cite{InfMatroidAxioms} for references. One of these, the `B-matroids' proposed by Higgs~\cite{Higgs69duality}, were later shown by Oxley~\cite{Oxley92} to describe%
   \COMMENT{}
   the models of any theory of infinite matroids that admitted both duality and minors as we know them. However, Higgs did not present his `B-matroids' in terms of axioms similar to those for finite matroids. As a consequence, theorems about finite matroids whose proofs rested on these axioms could not be readily extended to infinite matroids, even when this might have been possible in principle.

With the axioms from~\cite{InfMatroidAxioms} presented in the next section, this could now change: it should be possible now to extend many more results about finite matroids to infinite matroids, either by
\begin{itemize}
\item adapting their proofs based on the finite axioms to the (very similar) new infinite axioms,
\end{itemize}
or by
\begin{itemize}
\item finding a sequence of finite matroids that has the given infinite matroid as a limit, and is chosen in such a way that the instances of the theorem known for those finite matroids imply a corresponding assertion for the limit matroid.
\end{itemize}

After presenting our new axioms in Section~\ref{sec:axioms}, we apply them in Sections~\ref{sec:graphs} and~\ref{sec:proof} to see what they mean for graphs. We shall see that, for matroids whose circuits are the (usual finite) cycles of a graph, our axioms preserve what would be wrecked by the finitary axiom~(I4): that their duals are the matroids whose circuits are the bonds of our graph~-- even though these can now be infinite.

The converse is also nice. The dual $M^*$ of the (finitary) matroid $M$ whose circuits are the finite bonds of an infinite graph cannot be finitary; indeed, trivial exceptions aside, the duals of finitary matroids are never finitary~\cite{LasVergnas71,Bean76}; see also~\cite{BruhnWollanConInfMatroids}. So $M^*$ will have to have infinite circuits.%
   \footnote{We are using here that every dependent set contains a minimal such. This is indeed true.}%
   \COMMENT{}
   Excitingly, these circuits turn out to be familiar objects: when the graph is locally finite, they are the edge sets of the topological circles in its Freudenthal compactification, which already have a fixed place in infinite graph theory quite independently of matroids~\cite{RDsBanffSurvey}.

We shall see further that, for planar graphs, matroid duality is now fully compatible with graph duality as explored in~\cite{duality}. And Whitney's theorem, that a graph is planar if and only if its cycle matroid has a graphic dual, now has an infinite version too.

In some infinite graphs~$G$, including all locally finite ones, the elementary algebraic cycles%
   \footnote{These are the minimal 1-chains (possibly infinite) with zero boundary: the edge sets of finite cycles and of 2-way infinite paths.}
   form the circuits of a matroid, the {\em algebraic cycle matroid\/} of~$G$. We introduce this matroid, which was already studied by Higgs~\cite{Higgs69graphs}, in Section~\ref{sec:graphs}. In Section~\ref{sec:dualatc}, we determine its dual.

In Section~\ref{sec:thinsums}, finally, we consider representability. Infinite matroids that are representable in the usual sense are finitary, and hence cannot have representable duals. We suggest an adapted notion of representability based on infinite sums of functions to a field, and establish a sufficient condition for when this defines a matroid. Examples include the algebraic cycle matroids of graphs introduced in Section~\ref{sec:graphs}, and the algebraic cycle matroids of higher-dimensional complexes~\cite{InfMatroidAxioms}.

\section{Axioms}\label{sec:axioms}

We now present our five sets of axioms for finite or infinite matroids, in terms of independent sets, bases, circuits, closure and rank. These axioms were first stated in~\cite{InfMatroidAxioms}, and proved to be equivalent to each other in the usual sense. For finite or finitary matroids they coincide with the usual finite matroid axioms.

Let $E$ be any set, finite or infinite; it will be the default ground set for all matroids considered in this paper. We write $2^E$ for its power set. The set of all pairs $(A,B)$ such that $B\subseteq A\subseteq E$ will be denoted by $(2^E\times 2^E)_\subseteq$; for its elements we usually write $(A|B)$ instead of~$(A,B)$.%
   \COMMENT{}
   Unless otherwise mentioned, the terms `minimal' and `maximal' refer to set inclusion. Given $\E\sub 2^E$, we write $\Emax$ for the set of maximal elements of~$\E$, and $\dcl{\E}$ for the {\em down-closure} of~$\E$, the set of subsets of elements of~$\E$. For $F\sub E$ and $x\in E$, we abbreviate $F\sm\{x\}$ to $F-x$ and $F\cup\{x\}$ to~$F+x$. We shall not distinguish between infinite cardinalities and denote all these by~$\infty$; in particular, we shall write $|A|=|B|$ for any two infinite sets $A$ and~$B$. The set $\N$ contains~0.

One central axiom that features in all our axiom systems is that every independent set extends to a maximal one, even inside any restriction $X\sub E$.%
   \footnote{Interestingly, we shall not need to require that every dependent set contains a minimal dependent set. We need that too, but it follows~\cite{InfMatroidAxioms}.}
    The notion of what constitutes an independent set, however, will depend on the type of axioms under consideration. We therefore state this extension axiom in a more general form first, without reference to independence, so as to be able to refer to it later from within different contexts.

Let $\I\sub 2^E$. The following statement describes a possible property of~$\I$.

\begin{itemize}
\item[\rm (M)] Whenever $I\sub X\sub E$ and $I\in\I$, the set $\{\,I'\in\I\mid I\sub I'\sub X\,\}$ has a maximal element.
\end{itemize}
Note that the maximal superset of $I$ in~$\I\cap 2^X$ whose existence is asserted in~(M) need not lie in~$\Imax$.

\subsection{Independence axioms}\label{independenceaxioms}

The following statements about a subset $\I$ of $2^E$ are our {\em independence axioms\/}:
\begin{itemize}
\item[\rm (I1)] $\emptyset\in\mathcal I$.
\item[\rm (I2)] $\dcl{\I}=\I$, that is, $\I$ is closed under taking subsets.
\item[\rm (I3)] For all $I\in\I\sm\Imax$ and $I'\in\Imax$, there is an $x\in I'\sm I$ such that $I+x\in\I$.\looseness=-1
\item[\rm (IM)] $\I$ satisfies~(M).
\end{itemize}

\goodbreak

When a set $\I\sub 2^E$ satisfies the independence axioms, we call the pair $(E,\I)$ a {\em matroid\/} on~$E$. We call every element of $\I$ an {\em independent\/} set, every element of $2^E\sm\I$ a {\em dependent\/} set, the maximal independent sets {\em bases}, and the minimal dependent sets {\em circuits}. This matroid is {\em finitary\/} if it also satisfies (I4) from the Introduction, which is equivalent to requiring that every circuit be finite~\cite{InfMatroidAxioms}.

The $2^E\to 2^E$ function mapping a set $X\sub E$ to the set
 $$\cl(X) := X\cup \{\,x\mid \exists\, I\sub X\colon I\in\I\ \text{but}\ I+x\notin\I\,\} $$
will be called the {\em closure operator\/} on $2^E$ {\em associated with~$\I$.}

The $(2^E\times 2^E)_\subseteq\to \N\cup\{\infty\}$ function $r$ that maps a pair $A\supe B$ of subsets of $E$ to
 $$r(A|B) := \max\, \{\,|I\sm J| : I\supe J,\ I\in\I\cap 2^A,\ J\text{ maximal in }\I\cap 2^B\}$$
will be called the {\em relative rank function\/} on the subsets of $E$ {\em associated with~$\I$.} This maximum is always attained, and independent of the choice of~$J$~\cite{InfMatroidAxioms}.

\subsection{Basis axioms}

The following statements about a set $\B\sub 2^E$ are our {\em basis axioms\/}:
\begin{itemize}
\item[\rm (B1)] $\B\ne\es$.
\item[\rm (B2)] Whenever $B_1,B_2\in\B$ and $x\in B_1\sm B_2$, there is an element $y$ of~$B_2\sm B_1$ such that $(B_1-x)+y\in\B$.\COMMENT{}
\item[\rm (BM)] The set $\I:=\dcl{\B}$ of all {\em $\B$-independent\/} sets satisfies~(M).
\end{itemize}

\subsection{Closure axioms}

The following statements about a function $\cl\colon 2^E\to 2^E$ are our {\em closure axioms\/}:
\begin{itemize}
\item[\rm (CL1)] For all $X\sub E$ we have $X\sub \cl(X)$.
\item[\rm (CL2)] For all $X\sub Y\sub E$ we have $\cl(X)\sub \cl(Y)$.
\item[\rm (CL3)] For all $X\sub E$ we have $\cl(\cl(X)) = \cl(X)$.
\item[\rm (CL4)] For all $Z\sub E$ and $x,y\in E$, if $y\in\cl(Z+x)\sm\cl(Z)$ then $x\in\cl(Z+y)$.
\item[\rm (CLM)] The set $\I$ of all \cl-{\em independent\/} sets satisfies~(M). These are the sets $I\sub E$ such that $x\notin\cl(I-x)$ for all $x\in I$.
\end{itemize}

\subsection{Circuit axioms}

The following statements about a set $\C\sub 2^E$ are our {\em circuit axioms\/}:
\begin{itemize}
\item[\rm (C1)] $\es\notin\C$.
\item[\rm (C2)] No element of $\C$ is a subset of another.
\item[\rm (C3)] Whenever $X\sub C\in\C$ and $(C_x\mid x\in X)$ is a family of elements of~$\C$ such that $x\in C_y\Leftrightarrow x=y$ for all $x,y\in X$, then for every $z \in C\sm \left( \bigcup_{x \in X} C_x\right)$ there exists an element $C'\in\C$ such that $z\in C'\sub \left(C\cup  \bigcup_{x \in X} C_x\right) \sm X$.%
   \COMMENT{}
\item[\rm (CM)] The set $\I$ of all {\em $\C$-independent\/} sets satisfies~(M). These are the sets $I\sub E$ such that $C\not\sub I$ for all $C\in\C$.
\end{itemize}
Axiom (C3) coincides for $|X|=1$ with the usual (`strong')%
   \COMMENT{}
   circuit elimination axiom for finite matroids. In particular, it implies that adding an element to a basis creates at most one circuit;%
   \COMMENT{}
   the fact that it does create such a ({\em fundamental\/}) circuit is trivial when bases are defined from these circuit axioms (as maximal sets not containing a circuit), while if we start from the independence axioms it follows from the fact, mentioned before, that every dependent set contains a minimal one~\cite{InfMatroidAxioms}. We remark that the usual finite circuit elimination axiom is too weak to guarantee a matroid~\cite{InfMatroidAxioms}.

\subsection{Rank axioms}

The following statements about a function $r\colon (2^E\times 2^E)_\subseteq\to \N\cup\{\infty\}$ are our (relative) {\em rank axioms\/}:
\begin{itemize}
\item[\rm (R1)] For all $B\sub A\sub E$ we have $r(A|B)\le |A\sm B|$.
\item[\rm (R2)] For all $A,B\sub E$ we have $r(A|A\cap B)\ge r(A\cup B|B)$.
\item[\rm (R3)] For all $C\sub B\sub A\sub E$ we have $r(A|C) = r(A|B) + r(B|C)$.
\item[\rm (R4)] For all families $(A_\gamma)$ and $B$ such that $B\sub A_\gamma\sub E$ and $r(A_\gamma|B)=0$ for all~$\gamma$, we have $r(A|B)=0$ for $A := \bigcup_\gamma A_\gamma$.
\item[\rm (RM)] The set $\I$ of all $r$-{\em independent\/} sets satisfies~(M). These are the sets $I\sub E$ such that $r(I|I-x) > 0$ for all $x\in I$.
\end{itemize}

\goodbreak

For finite matroids, these axioms (with (R4) and (RM) becoming redundant) are easily seen to be tantamount to the usual axioms for an absolute rank function~$R$ derived as $R(A) := r(A|\es)$, or conversely with ${r(A|B) := R(A) - R(B)}$%
   \COMMENT{}
   for $B\sub A$.%
   \COMMENT{}

\section{Bond and cycle matroids}\label{sec:graphs}

In this section we develop the theory of our axioms to see what it yields for the usual matroids for graphs when these are infinite. See~\cite{InfMatroidAxioms} for applications to other structures than graphs. All our graphs may have parallel edges and loops.

A well-known matroid associated with a finite
graph $G$ is its cycle matroid: the matroid whose circuits are the edge sets of the cycles in~$G$. The bases of this matroid are the edge
sets of the {\em spanning forests\/} of~$G$, the sets that form a spanning tree in every component of~$G$. This construction works in infinite 
graphs too: the edge sets of the finite cycles in $G$ form the circuits of a finitary matroid~$\mfc(G)$, whose bases are the edge sets of the spanning forests of~$G$. We shall call $\mfc(G)$ the \emph{finite-cycle matroid of $G$}. Similarly, we let the 
\emph{finite-bond matroid} $\mfb(G)$ of $G$ be the matroid whose circuits are the finite bonds of~$G$. (A~{\em  bond\/} is a minimal non-empty cut.)
This, too, is a finitary matroid. 

If $G$ is finite, then $\mfc(G)$ and $\mfb(G)$ are dual to each other. For infinite~$G$, however, things are different. As remarked earlier, the duals of finitary matroids are not normally finitary, so the duals of $\mfc(G)$ and $\mfb(G)$ will in general have infinite circuits. In the case of~$\mfc(G)$, its cocircuits are the expected ones, the (finite or infinite) bonds:

\begin{theorem}\label{easythm}
Let $G$ be any graph.
\begin{enumerate}[\rm (i)]
\item The bonds of $G$, finite or infinite, are the circuits of a matroid~$\mb(G)$.
\item This matroid is the dual of the finite-cycle matroid~$\mfc(G)$ of~$G$.
\end{enumerate}
\end{theorem}

\noindent
The matroid $\mb(G)$ defined in Theorem~\ref{easythm} will be called the {\em bond matroid} of~$G$. We defer the proof of the theorem to Section~\ref{sec:proof}; it is essentially the same as for finite graphs, although now the bonds can be infinite.

\bigskip

Similarly, the dual of $\mfb(G)$ will in general have infinite circuits. Ideally, these would form some sort of `infinite cycles' in~$G$. `Infinite cycles' have indeed been considered before for graphs, though in a purely graph-theoretic context: there is a topological such notion that makes it possible to extend classical results about cycles in finite graphs (such as Hamilton cycles) to infinite graphs, see~\cite{RDsBanffSurvey} and~\cite{MSsBanffSurvey} in this issue. Rather strikingly, it turns out that these `infinite cycles' are the solution also to our problem: their edge sets are precisely the (possibly infinite) cocircuits of~$\mfb(G)$.

In order to define those `infinite cycles', we need to endow our given graph~$G$
with a topology. A \emph{ray} is a one-way infinite path. Two rays are \emph{edge-equivalent}
if for any finite set $F$ of edges there is a component of $G-F$ that contains
subrays of both rays. The equivalence classes of this relation are 
the \emph{edge-ends of $G$}; we denote the set of these edge-ends by~$\eends(G)$.

Let us view the edges of~$G$ as disjoint topological copies of~$[0,1]$, and let $X_G$ be the quotient space obtained by identifying these copies at their common vertices. The set of inner points of an edge~$e$ will be denoted by~$\kreis e$. We now define a topological space $\etop G$ on the point set of $X_G\cup\eends(G)$ by taking as our open sets the unions of sets $\widetilde C$, where $C$ is a connected component of $X_G - Z$ for some finite set $Z\subset X_G$ of inner points of edges, and $\widetilde C$ is obtained from $C$ by adding all the edge-ends represented by a ray in~$C$.%
   \footnote{Note that $\etop G$ induces on $G$ a topology coarser than that of~$X_G$. We needed $X_G$ only provisionally, so that the connected components used in the definition of~$\etop G$ were defined.}

When $G$ is connected, $\etop G$~is a compact topological space~\cite{Moritz}, although in general it need not be Hausdorff: the common starting vertex of infinitely many otherwise disjoint equivalent rays, for example, cannot be distinguished topologically from the edge-end which those rays represent, and neither can two vertices joined by infinitely many edges or edge-disjoint paths.%
   \footnote{A consequence that takes some getting used to is that circles in $\etop G$ need not look `round': an edge joining two indistinguishable vertices, for example, will form a topological loop with either one of them.}
   However if $G$ is locally finite, then $\etop{G}$ coincides with the (Hausdorff) Freudenthal compactification of~$G$. See Section~\ref{sec:proof} for more properties of~$\etop G$.

For any set $X\subseteq\etop{G}$ we call
 $$E(X):=\{e\in E(G):\kreis e\subseteq X\}$$
the \emph{edge set of $X$}.
A subspace $C$ of $\etop G$ that is homeomorphic to~$S^1$
is a 
 {\em circle\/} in~$\etop G$.
   % grrr, circles are not standard subspaces, as they are not necessarily closed...
One can show that $\bigcup E(C)$ is dense in~$C$, so $C$ lies in%
   \COMMENT{}
    the closure of the subgraph formed by its edges~\cite{Moritz}. In particular, there are no circles consisting only of edge-ends.

A subspace $X\sub\etop G$ is a {\em standard subspace} if it is the closure in $\etop{G}$ of a subgraph of~$G$.
A \emph{topological spanning tree} of $G$ is a standard subspace $T$ of $\etop{G}$ that is path-connected and contains $V(G)$ but contains no circle. Note that, since standard subspaces are closed, $T$~will also contain~$\eends(G)$.

\goodbreak

\begin{theorem}\label{missingdual}
Let $G$ be any connected\/%
   \footnote{The theorem extends to disconnected graphs in the obvious way.}\!%
   \COMMENT{}
   graph.
\begin{enumerate}[\rm (i)]
\item The edge sets of the circles in~$\etop G$ are the circuits of a matroid~$\mc(G)$, the {\em cycle matroid} of~$G$.
\item The bases of $\mc(G)$ are the edge sets of the topological spanning trees of~$G$.
\item The cycle matroid~$\mc(G)$ is the dual of the finite-bond matroid~$\mfb(G)$.
\end{enumerate}
\end{theorem}

\noindent
We shall prove Theorem~\ref{missingdual} in Section~\ref{sec:proof}.

\medbreak

In the finite world, matroid duality is compatible with graph duality in that the dual of the cycle matroid of a finite planar graph $G$ is the cycle matroid of its (geometric or algebraic) dual~$G^*$. Duality for infinite graphs has come to be properly understood only recently~\cite{duality}. But now that we have matroid duality as well, it turns out that the two are again compatible. In the remainder of this section we briefly explain how infinite graph duality is defined, and then show its compatiblity with matroid duality.

When one tries to define abstract graph duality so that it satisfies the minimum requirement of capturing the geometric duality of locally finite graphs in the plane (where one has a dual vertex for every face and a dual edge between vertices representing two faces for every edge that lies on the boundary of both these faces), the first thing one realizes is that by taking duals one will leave the class of locally finite graphs: the dual of a ray, for example, is a vertex with infinitely many loops. On the other hand, Thomassen~\cite{thomassen82} showed that any class of graphs for which duality can be reasonably defined cannot be much larger: these graphs have to be \emph{finitely separable} in that every two vertices can be separated by finitely many edges.%
   \footnote{Christian, Richter and Rooney~\cite{CRR09} define certain dual objects for arbitrary planar graphs; however these objects are `graph-like spaces', not graphs.}

It was finally shown in~\cite{duality} that the class of finitely separable graphs is indeed the right setting for infinite graph duality, defined as follows. Let $G$ be a finitely separable graph. A graph $G^*$ is called a \emph{dual of $G$} if there is a bijection
 $${}^*:E(G)\to E(G^*)$$
such that a set $F\subseteq  E(G)$ is the edge set of a circle in $\etop{G}$ if and only
if $F^* := \,\{e^*\mid e\in F\/\}$ is a bond of $G^*$.%
   \footnote{We are cheating a bit here, but only slightly. In~\cite{duality}, these circles are taken not in~$\etop{G}$ but in a slightly different space $\tilde G$. However, while the circles in $\tilde G$ may differ slightly from those in~$\etop{G}$, their edge sets are the same. This is not hard to see directly; it also follows from Theorems 6.3 and~6.5 in~\cite{tst} in conjunction with Satz 4.3 and~4.5 in~\cite{Moritz}.}%
   \COMMENT{}
   Duals defined in this way behave just as for finite graphs:

\begin{theorem}\selfcite{duality}\label{graphduality}
Let $G$ be a countable finitely separable graph.
\begin{enumerate}[\rm (i)]
\item  $G$ has a dual if and only if $G$ is planar.
\item If $G^*$ is a dual of $G$, then $G^*$ is finitely separable, $G$ is a dual of $G^*$, and this is
witnessed by the inverse bijection of~${}^*$.
\item Duals of 3-connected graphs are unique, up to isomorphism.
\end{enumerate}
\end{theorem}

\goodbreak

At the time, the reason for defining graph duality as above was purely graph-theoretic: it appeared (and still appears) to be the unique  way to make all three statements of Theorem~\ref{graphduality} true for infinite graphs. As matroid duality was developed independently of graph duality, it is thus remarkable---and adds to the justification of both notions---that the two are once more compatible, as far as remains possible in an infinite setup:

\begin{theorem}
\label{graphdual}
Let $G$ and $G^*$ be a pair of countable  dual graphs, each finitely separable, and defined on the same edge set~$E$. Then 
\[
\mfb(G)=\mc^*(G)=\mb^*(G^*)=\mfc(G^*).
\]
\end{theorem}

\begin{proof}
The first equality is Theorem~\ref{missingdual}\,(iii). The last equality is Theorem~\ref{easythm}\,(ii) (after dualizing). The middle equality follows from $\mc(G)=\mb(G^*)$, which is a direct consequence of the definition of a dual graph.
\end{proof}

Finally, we obtain an infinite analogue of Whitney's theorem that a finite graph is planar if and only if the dual of its cycle matroid is `graphic', i.e., is the cycle matroid of another finite graph. In our more general context, let us call a matroid \emph{graphic} if it is isomorphic to the cycle matroid of some graph, and \emph{finitely graphic} if it is 
isomorphic to the finite-cycle matroid of some graph.

\goodbreak

\begin{theorem}\label{graphicduals}
The following three assertions are equivalent for a countable finitely separable graph~$G$:
\begin{enumerate}[\rm (i)]
\item $G$ is planar;
\item $\mc^*(G)$ is finitely graphic;
\item $\mfc^*(G)$ is graphic.
\end{enumerate}
\end{theorem}

\noindent
We shall prove Theorem~\ref{graphicduals} in Section~\ref{sec:proof}. The equivalence of (i) and~(ii) can also be derived from a more general result of Christian et al.\ on `graph-like spaces'~\cite{CRR09}.

\medbreak

Another natural matroid in a locally finite graph $G$ is its algebraic cycle matroid: the matroid whose circuits are the {\em elementary algebraic cycles\/} of $G$, the minimal non-empty edge sets inducing even degrees at all the vertices. Clearly, these are the edge sets of the finite cycles in $G$ and those of its {\em double rays\/}, its 2-way infinite paths.

The elementary algebraic cycles do not form a matroid in every infinite graph: it is easy to check \cite[Section~6]{InfMatroidAxioms} that they do not satisfy our circuit axioms when $G$ is the {\em Bean graph\/} shown in Figure~\ref{Beangraph}. However, by a result of Higgs~\cite{Higgs69graphs} made applicable to our matroids by \cite[Theorem 5.1]{InfMatroidAxioms}, this is essentially the only counterexample:

\begin{figure}[htbp]
  \centering
%\epsfxsize=0.7\hsize  
%\noindent
%\epsfbox{Beangraph.eps}
  \includegraphics[width=0.7\linewidth]{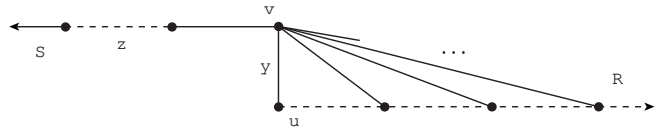}
  \caption{The Bean graph}
  \label{Beangraph}
\end{figure}

\begin{theorem}[Higgs 1969]\label{HiggsBean}
The elementary algebraic cycles of an infinite graph $G$ are the circuits of a matroid on its edge set $E(G)$ if and only if $G$ contains no subdivision of the Bean graph.
\end{theorem}

\begin{corollary}
The elementary algebraic cycles of any locally finite graph are the circuits of a matroid.\qed
\end{corollary}

We call the matroid from Theorem~\ref{HiggsBean} the {\em algebraic cycle matroid\/}~$\atc(G)$ of the graph~$G$. In Section~\ref{sec:dualatc} we determine its dual~$\atc^*(G)$: it is the matroid whose circuits are the minimal non-empty cuts of $G$ at least one side of which contains no ray.%
   \COMMENT{}

The algebraic cycle matroid is {\it representable\/} in a sense adapted to non-finitary matroids, which we discuss in Section~\ref{sec:thinsums}. We shall also prove a general sufficient condition for this notion of representability.

\goodbreak

\section{Proofs of Theorems \ref{easythm}, \ref{missingdual} and~\ref{graphicduals}}\label{sec:proof}

We begin with the easy proof of Theorem~\ref{easythm}, which we restate:

\proclaim Theorem~\ref{easythm}. Let $G$ be any graph.
\begin{enumerate}[\rm (i)]
\item The bonds of $G$, finite or infinite, are the circuits of a matroid $\mb(G)$, the {\em bond matroid} of~$G$.
\item The bond matroid of~$G$ is the dual of its finite-cycle matroid~$\mfc(G)$.
\end{enumerate}

\vskip-\medskipamount\vskip0pt
\begin{proof}
For simplicity we assume that $G$ is connected; the general case is very similar. From~\cite{InfMatroidAxioms} we know that $\mfc(G)$ has a dual; let us call this dual~$\mb(G)$, and show that its circuits are the bonds of~$G$. By definition of matroid duality, the circuits of $\mb(G)$ are the minimal edges sets that meet every spanning tree of~$G$.

We show first that every bond $B$ of $G$ is a circuit of~$\mb(G)$, a minimal set of edges meeting every spanning tree. Since $B$ is a non-empty cut, it is the set of edges across some partition of the vertex set of~$G$. Every spanning tree meets both sides of this partition, so it has an edge in~$B$. On the other hand, we can extend any edge $e\in B$ to a spanning tree of $G$ that contains no further from~$B$, since by the minimality of $B$ as a cut its two sides are connected in~$G$. Hence $B$ is minimal with the property of meeting every spanning tree.

Conversely, let $B$ be any set of edges that is minimal with the property of meeting every spanning tree. We show that $B$ contains a bond; by the implication already shown, and its minimality, it will then {\em be\/} that bond. Since $G$ has a spanning tree, we have $B\ne\emptyset$; let $e\in B$. If $B$ contains no bond, then every bond has an edge not in~$B$. The subgraph $H$ formed by all these edges is connected and spanning in~$G$, as otherwise the edges of $G$ from the component $C$ of $H$ containing~$e$ to any fixed component of $G-C$ would form a bond of $G$ with no edge in~$H$, contradicting its definition. So $H$ contains a spanning tree. This misses~$B$, contradicting the choice of~$B$.
\end{proof}

We prove Theorem~\ref{missingdual} for countable graphs; the proof for arbitrary graphs can be deduced from this by considering a quotient space of $\etop G$ as explained in~\cite{Moritz}.
   \COMMENT{}
   For the remainder of this section, let $G$ be a fixed countable%
   \COMMENT{}
   connected graph.\looseness=-1

We shall call two points in $\etop{G}$ \emph{(topologically) indistinguishable} if they have the same open neighbourhoods. Clearly two vertices or edge-ends $x,y\in\etop{G}$ are indistinguishable if they cannot be separated by finitely
many edges. (If both are edge-ends, then $x=y$.) On the other hand, two such points that can be separated by finitely 
many edges have disjoint open neighbourhoods. Inner points of edges are always distinguishable from all other points.

We shall need a few lemmas. Some of these are quoted from Schulz~\cite{Moritz}; the others are adaptations of results proved in~\cite{tst} for the special case that $G$ is finitely separable. We remark that it is also possible to reduce Theorem~\ref{missingdual} formally to that case by replacing $G$ with a quotient graph as explained in~\cite{Moritz}.

\begin{lemma}{\bf \cite{Moritz}}\label{etopcomp}
$\etop{G}$~is a compact space. 
\end{lemma}

\begin{lemma}\label{opensets}
Let $X\subseteq\etop{G}$ be a closed subspace, with disjoint open subsets $O_1,O_2$ such that $X=O_1\cup O_2$.%
   \COMMENT{}
   Then the set $F$ of edges of $G$ with one endvertex in~$O_1$ and the other in~$O_2$ is finite.
\end{lemma}
\begin{proof}
Suppose $F$ is infinite. As $\etop G$ is compact, the mid-points of the edges in~$F$ have an accumulation point $x$ in~$\etop G$. By definition of~$\etop G$, every neighbourhood of $x$ contains infinitely many edges from~$F$,%
   \COMMENT{}
   and hence meets both $O_1$ and~$O_2$. Since $O_1$ and $O_2$ are closed in~$X$, and hence in~$\etop G$, this means that $x\in O_1\cap O_2$, contradicting our assumption that $O_1\cap O_2 = \emptyset$.
 \end{proof}

In a Hausdorff space, every topological $x$--$y$ path contains an injective such path, an $x$--$y$ {\em arc\/}. Since $\etop G$ is not necessarily Hausdorff we cannot assume this shortcut lemma in general, but it holds in the relevant case:

\begin{lemma}{\bf \cite{Moritz}}\label{patharc}
If two points $x,y\in V(G)\cup\eends(G)$ are separated by a finite set of edges, then every
topological $x$--$y$ path contains an $x$--$y$ arc.
\end{lemma}

\begin{lemma}{\bf \cite{Moritz}}\label{pathlimit}
Let $x,y\in V(G)\cup\eends(G)$,
and let $(A_\gamma)_{\gamma<\lambda}$ be a transfinite
sequence of $x$--$y$ arcs in $\etop{G}$. Then there exists a topological $x$--$y$ path $P$
and a dense subset $P^*$ of $P$ so that for all $p\in P^*$ the arcs
$A_\gamma$ containing $p$ form a cofinal subsequence.
\end{lemma}

\begin{lemma}\label{pathcon}
Every closed connected subspace $X$ of $\etop{G}$ is path-connected.
\end{lemma}

%This is basically the proof of Diestel and K\"uhn~\cite{tst}.
\begin{proof}
Suppose $X$ is connected but not path-connected. Then there are $x,y\in V(G)\cup\eends(G)$
contained in different path-components. In particular, $x$ and $y$ are topologically distinguishable,%
   \COMMENT{}
   so they are separated by finitely many edges.
Let $e_1,e_2,\ldots$ be a (possibly finite) enumeration of the edges in $E(G)\sm E(X)$,
let $F_i:=\{e_1,\ldots,e_i\}$ for all $i$. If there exists an $i$ such that $x$ and $y$ lie in the closures of different graph-theoretical components of $G-F_i$, then picking an inner point outside $X$ from every edge in $F_i$ we obtain a finite set $Z\subseteq\etop G\sm X$ witnessing that $x$ and $y$ lie in distinct open sets of $X$ whose union is all of~$X$,%
   \COMMENT{}
   contradicting our assumption that $X$ is connected.

Hence for every~$i$ the points $x$ and $y$ lie in the closure $\overline C_i$ of the same component $C_i$ of ${G - F_i}$. So for each $i$ there is a path, ray or double ray  connecting $x$ to $y$ in~$\overline C_i$,
and with Lemma~\ref{patharc} we then obtain an $x$--$y$ arc $A_i$ in~$\overline C_i$. 
By Lemma~\ref{pathlimit} this implies that there is a topological
$x$--$y$ path $P$ and a dense subset $P^*\subseteq P$ such that for every $p\in P^*$ 
the arcs $A_i$ containing $p$ form a cofinal subsequence. Suppose there exists a $j$
such that $\kreis e_j\subseteq P$. Then there must be a point $p\in \kreis e_j\cap P^*$.
However, none of the $A_i$ with $i\geq j$ contains~$\kreis e_j$. Thus, $P$ 
does not use any edge outside $X$. As $X$ is closed, this implies that $P\subseteq X$. The required $x$--$y$ arc in $X$ can be found inside~$P$ by Lemma~\ref{patharc}.
\end{proof}

\begin{lemma}\label{extension}
Let $F\subseteq E(G)$ be a set of edges whose closure in $\etop G$ contains no circle. Then $G$ has a topological spanning tree whose edge set contains~$F$.
\end{lemma}

\proof
Let $G=(V,E)$, let $e_1,e_2,\ldots$ be an enumeration of the edges in $E\sm F$, and set 
$T_0:=E$. Inductively, if the closure of $(V,T_{i-1}-e_i)$ is 
connected in $\etop{G}$ then set $T_i:=T_{i-1}-e_i$;
otherwise put $T_i:=T_{i-1}$. Finally, we set $T:=\bigcap_{i=0}^{\infty} T_i$.

In order to show that $T$ is the edge set of a topological spanning tree,
let us first check that the closure $X$ of $(V,T)$ is connected.
Suppose there are two disjoint non-empty open sets $O_1$ and $O_2$ of $X$ with 
$X=O_1\cup O_2$. Then Lemma~\ref{opensets} implies that the cut $S$ 
consisting of the edges with one endvertex in $O_1$ and the other in $O_2$
is finite. If $j$ is the largest integer with $e_j\in S$ then, however, 
the closure of $(V,T_j)$ is not connected, a contradiction.
Thus, $\overline T=X$ is connected and therefore spanning. Moreover, 
$\overline T$ is path-connected, by Lemma~\ref{pathcon}.

Secondly, we need to show that $\overline T$ is acirclic. So, suppose that 
$\overline T$ contains a circle $C$. Since every circle lies in the closure of its edges but the closure of $\bigcup F$ contains no circle, $E(C)\sm F$ is non-empty. 
Pick $j$ minimal with $e_j\in E(C)\sm F$. Since $e_j$ was not deleted from~$T_{j-1}$ when $T_j$ was formed, 
the closure $Y$ of $(V,T_{j-1}-e_j)$ is disconnected. So there are two disjoint non-empty
open subsets $O_1,O_2$ of $Y$ such that $Y=O_1\cup O_2$. The endvertices of $e_j$ do not lie in the same~$O_i$, since adding $e_j$ to that $O_i$ would then yield a similar decomposition of the closure of $(V,T_{j-1})$, contradicting its connectedness.
But now the connected subset $C\sm \kreis e_j$ of $Y$ meets both $O_1$ and $O_2$, 
a contradiction. Thus, $\overline T$ does not contain any circle
and is therefore a topological spanning tree.
\endproof

\begin{lemma}\label{incomp}
Let $C_1$ and $C_2$ be two circles in $\etop{G}$.
Then $E(C_2)\subseteq E(C_1)$ implies that $E(C_1)=E(C_2)$.
\end{lemma}
\begin{proof}

We first prove the following: 
\begin{equation}\label{boundary}
\begin{minipage}[c]{0.8\textwidth}\em
For every point $x\in\overline{C_1}\sm C_1$ there is a point $y\in C_1$ 
such that $x$ and $y$ are indistinguishable.
\end{minipage}\ignorespacesafterend 
\end{equation} 
Indeed, consider a $z\in\etop{G}$ that is distinguishable 
from all points in $C_1$. Thus, we may pick for every $p\in C_1$
two disjoint open neighbourhoods $O_z^p$ and $O_p$ of $z$ and $p$, respectively. 
Note that
$C_1$ is compact, being a continuous image of the compact space~$S^1$.
Thus, there is a finite subcover $O_{p_1}\cup\ldots\cup O_{p_n}$ of $C_1$.
Then, the open set $\cap_{i=1}^n O_z^{p_i}$ is disjoint from $C_1$ and contains $z$.
Hence, $z$ does not lie in the closure of $C_1$. This proves~\eqref{boundary}.

Next, suppose that $E(C_2)$ is a proper subset of $E(C_1)$, and pick $e\in E(C_2)$ and $f\in E(C_1)\sm E(C_2)$.
Since $X:=C_1\sm(\kreis e\cup\kreis f)$ is disconnected there exist two disjoint non-empty open sets $O'_1$ and $O'_2$
of $X$ with $X=O'_1\cup O'_2$. For $j=1,2$, denote by $I_j$ the set of points $x$ in $\etop{G}$ for which there is a 
 $y\in O'_j$ such that $x$ and $y$ are indistinguishable.
Then $O_1:=O'_1\cup I_1$ and $O_2:=O'_2\cup I_2$ are disjoint and open subsets of $X\cup I_1\cup I_2$.
Moreover, it follows from~\eqref{boundary} that $\overline{C_1}\sm(\kreis e\cup\kreis f)=X\cup I_1\cup I_2$.
Therefore,   $O_1$ and $O_2$ are two disjoint non-empty open sets
of $\overline{C_1}\sm(\kreis e\cup\kreis f)$
with $\overline{C_1}\sm(\kreis e\cup\kreis f)=O_1\cup O_2$. 

As $C_2\sm\kreis e$ is a connected subset of $\overline{C_1}\sm(\kreis e\cup\kreis f)$ it lies in $O_1$ or in $O_2$, let us say in $O_1$. Then $\tilde O_1:=O_1\cup\kreis e$ and $O_2$ are two disjoint non-empty open subsets of $\overline{C_1}\sm\kreis f$ with $\overline{C_1}\sm\kreis f=\tilde O_1\cup O_2$. By~\eqref{boundary}, this means that also $C_1\sm\kreis f$ is disconnected.%
   \COMMENT{}
   But $C_1\sm\kreis f$ is a continuous image of a connected space, and hence connected.
\end{proof}

\goodbreak

\begin{lemma}\label{minmaxtst}
Let $T$ be a standard subspace of  $\etop G$. 
Then the following statements are equivalent:
\begin{enumerate}[\rm (i)]
\item $T$ is a topological spanning tree of $\etop{G}$.
\item $T$ is maximally acirclic, that is, it does not contain a circle but 
adding any edge in $E(G)\sm E(T)$ creates one.
\item $E(T)$ meets every finite bond, and is minimal with this property.
\end{enumerate}
\end{lemma}

\begin{proof}
Let us first prove a part of (iii)$\rightarrow$(i) before dealing with
all the other implications.
\begin{equation}\label{clm2}
\begin{minipage}[c]{0.8\textwidth}\em
If $E(T)$ meets every finite bond then $T$ is spanning and path-connected.
\end{minipage}\ignorespacesafterend 
\end{equation} 
Suppose that the closure $X$ of $(V(G),E(T))$ is not connected. 
Then there are two disjoint non-empty open sets $O_1$ and $O_2$ of $X$
with $X=O_1\cup O_2$. From Lemma~\ref{opensets} we get that the 
cut consisting of the edges with one endvertex in $O_1$ and the other in $O_2$
is finite.
Since each of $O_1$ and $O_2$ needs to contain a vertex, this cut is
non-empty. Hence, $E(T)$ 
misses a finite bond, a contradiction. Therefore, $T=X$ is connected and then, by Lemma~\ref{pathcon},
path-connected. 

(i) $\rightarrow$ (ii) Consider any edge $e\notin E(G)\sm E(T)$. 
If the endvertices $u$ and $v$ of $e$ cannot be separated by finitely many
edges then $e-u$ (and also $e-v$) is a circle in $\etop{G}$. Otherwise, any topological
$u$--$v$ path contains an $u$--$v$ arc by Lemma~\ref{patharc}. 
In particular, $T$ contains an $u$--$v$ arc that together with $e$ forms a circle.

(ii) $\rightarrow$ (iii) Suppose that $E(T)$ misses a finite bond $F$.
Pick $e\in F$, and let $C$ be a circle in $T\cup e$ through $\kreis e$. 
Pick an inner point of every edge in $F$ and denote the set of these
points by $Z$. Then the two components of $\etop{G}\sm Z$, each of which
contains an endvertex of $e$, form two 
 disjoint open sets containing $T$.
However, $C\sm\kreis e\subseteq T$ is a connected set that meets both of these disjoint open sets,
which is impossible. 
Thus, $E(T)$ meets every finite
bond. In particular, $T$ is spanning and path-connected, by~\eqref{clm2}.

Let $f$ be any edge in $E(T)$, and let us show that $E(T)-f$ misses
some finite bond.
Denote
the endvertices of $f$ by $r$ and $s$, and observe that $r$ and $s$ 
can be separated by finitely many edges as $T$ is acirclic. 
Denote by $K_r$ and $K_s$
the path-components of $T\sm\kreis f$ containing $r$ and $s$, respectively. 
By Lemma~\ref{patharc} and 
as $T$ does not contain any circle, $K_r$ and $K_s$ are distinct, and thus
disjoint. As $T$ is path-connected, it follows that $T\sm\kreis f$
is the disjoint union of the open sets $K_r$ and $K_s$.
Now Lemma~\ref{opensets} yields that there are only finitely many edges with
one endvertex in $K_r$ and the other in $K_s$. As $T$ is spanning this
means that $E(T)-f$ misses a finite cut.

(iii) $\rightarrow$ (i) 
By~\eqref{clm2}, we only need to
check that $T$ does not contain any circle. Suppose there exists a circle $C\subseteq T$,
and pick some $e\in E(C)$. By the minimality of $E(T)$ there exists a finite bond $F$ so that
$F$ is disjoint from $T\sm\kreis e$. Then, however, picking inner points from the edges in 
$F$ yields a set $Z$, so that the connected set $C\sm\kreis e$ is contained
in $\etop{G}\sm Z$ but meets two components of $\etop{G}\sm Z$, 
which is impossible.
\end{proof}

We can now prove our main theorem, which we restate:

\proclaim Theorem~\ref{missingdual}.
\begin{enumerate}[\rm (i)]
\item The edge sets of the circles in~$\etop G$ are the circuits of a matroid~$\mc(G)$, the {\em cycle matroid} of~$G$.
\item The bases of $\mc(G)$ are the edge sets of the topological spanning trees of~$G$.
\item The cycle matroid~$\mc(G)$ is the dual of the finite-bond matroid~$\mfb(G)$.
\end{enumerate}

\vskip-\medskipamount\vskip0pt\begin{proof}
To bypass the need to verify any matroid axioms, we {\em define\/} $\mc(G)$ as the dual of~$\mfb(G)$ (which we know exists~\cite{InfMatroidAxioms}), ie., as the matroid whose bases $B$ are the complements of the bases of~$\mfb(G)$. These latter are the maximal edge sets not containing a finite bond, so the bases $B$ of $\mc(G)$ are the minimal edge sets meeting every finite bond. By Lemma~\ref{minmaxtst} below, this is equivalent to $B$ being the edge set of a topological spanning tree of~$\etop{G}$.%
   \COMMENT{}

We have defined $\mc(G)$ so as to make (iii) true, and shown~(ii). It remains to show~(i): that the circuits of $\mc(G)$ are the edge sets of the circles in~$\etop G$. Since no circuit of a matroid contains another circuit, and since by Lemma~\ref{incomp} no edge set of a circle contains another such set, it suffices to show that every circuit contains the edge set of a circle, and conversely every edge set of a circle contains a circuit.

For the first of these statements note that, by assertion~(ii), a circuit $D$ of~$\mc(G)$ does not extend to the edge set of a topological spanning tree.%
   \COMMENT{}
   Hence by Lemma~\ref{extension} its closure $\overline{\bigcup D}$ in $\etop G$ contains a circle~$C$. For the second statement, note that the edge set $D$ of a circle $C$ is not contained in the edge set of a topological spanning tree~$T$, because $T$ is closed and would therefore contain $\overline{\bigcup D}\supseteq C$, contradicting its definition. Hence $D$ is dependent in~$\mc(G)$, and therefore contains a circuit~\cite{InfMatroidAxioms}.
\end{proof}

Finally, let us restate and prove Theorem~\ref{graphicduals}:

\proclaim Theorem~\ref{graphicduals}.
The following three assertions are equivalent for a countable finitely separable graph~$G$:
\begin{enumerate}[\rm (i)]
\item $G$ is planar;
\item $\mc^*(G)$ is finitely graphic;
\item $\mfc^*(G)$ is graphic.
\end{enumerate}

\begin{proof}
If $G$ is planar, it has a dual~$G^*$. Then $\mc^*(G) =\mfc(G^*)$ by Theorem~\ref{graphdual}, and $\mfc^*(G) = \mb(G) = \mc(G^*)$ by Theorems \ref{easythm} and~\ref{graphduality}\,(ii).

(ii)$\to$(i): Since $\mc^*(G)$ is finitely graphic, there exists a graph $H$ with the same edge set as $G$
such that $\mc^*(G) =\mfc(H)$. As $\mfc^*(H) =\mb(H)$ by Theorem~\ref{easythm} and matroid duals are unique, we obtain
$\mc(G)=\mb(H)$. Hence the edge sets of the circles in~$\etop{G}$, which by Theorem~\ref{missingdual}\,(i) are the circuits of~$\mc(G)$, are precisely the bonds of~$H$. So $H$ is a dual of $G$, and $G$ is planar by Theorem~\ref{graphduality}\,(i).

(iii)$\to$(i): Let $H$ be a graph such that $\mc(H) = \mfc^*(G)$. To show that $G$ is planar, it suffices by Kuratowski's theorem%
   \footnote{Its extension to countable graphs is straightforward by compactness; see~\cite[Exercise~8.23]{DiestelBook10noEE}.}
   to check that $G$ has no $K_5$- or $K_{3,3}$-minor, or in matroid terms, that $\mfc(G)$ has no minor isomorphic to $\mc(K_5)$ or $\mc(K_{3,3})$. As $\mfc^*(G) = \mc(H)$, this is equivalent to saying that $\mc(H)$ has no $\mc^*(K_5)$ or $\mc^*(K_{3,3})$-minor.%
   \COMMENT{}
   These latter two matroids are not graphic \cite[Prop.~2.3.3]{OxleyBook}, so it thus suffices to show that the finite minors of $\mc(H)$ are graphic.%
   \COMMENT{}

To prove this, consider a finite minor $M$ of~$\mc(H)$, obtained by deleting the set $X\sub E(H)$ and contracting the set $Y\sub E(H)$, say. Let $V$ be the finite set of vertices of $H$ incident with an edge in the ground set $E$ of~$M$, and let $K$ be the finite graph obtained from the graph $(V,E)$ by identifying any two vertices that are either indistinguishable in~$\etop H$ or joined by an arc in $\etop H$ whose edges lie in some fixed base $B$ of the restriction of $\mc(H)$ to~$Y$. Using Lemma~\ref{patharc} and \cite[Lemma~3.5]{InfMatroidAxioms}, it is now easy to show that $M = \mc(K)$. Hence $M$ is graphic, as desired. 
\end{proof}

We remark that the graphs witnessing (ii) and (iii) in Theorem~\ref{graphicduals} can be chosen to be finitely separable, too. Indeed, the graph $G^*$ which we used in our proof as a witness for both (ii) and (iii) is finitely separable by Theorem~\ref{graphduality}\,(ii).

\section{The dual of the algebraic cycle matroid}\label{sec:dualatc}

Recall from Theorem~\ref{HiggsBean} that the elementary algebraic cycles of a graph $G$ are the circuits of a matroid, the {\em algebraic cycle matroid\/}~$\atc(G)$ of~$G$, if and only if $G$ contains no subdivision of the graph shown in Figure~\ref{Beangraph}. In this section we characterize their matroid~$\atc^*(G)$.\looseness=-1

Recall that a \emph{ray} is a $1$-way infinite path. Let us call a non-empty cut $F=E(A,B)$ of $G$ \emph{\xsided}  if one of its sides $A,B$ is {\em small\/} in the sense that the subgraph it induces in $G$ contains no ray and $F$ is minimal with this property among the non-empty cuts of~$G$. If $G$ is connected then so is the small side of any \xsided\ cut, so this will be finite if $G$ is connected and locally finite.

\begin{theorem}\label{finsidedthm}
The cocircuits of a matroid $M$ that is the algebraic cycle matroid of a graph~$G$ are precisely the \xsided\ cuts of~$G$.
\end{theorem}

Casteels and Richter~\cite{CR08} studied a related duality: they showed that, in a locally finite
graph, the cuts with a finite side form the orthogonal space of the 
set of elementary algebraic cycles.%
   \COMMENT{}

For our proof of Theorem~\ref{finsidedthm} we need the following easy lemma from~\cite{InfMatroidAxioms}:

\begin{lemma}\label{circuitcocircuitcap}
A circuit and a cocircuit of a matroid never meet in exactly one element.
\end{lemma}

\begin{proof}[Proof of Theorem~\ref{finsidedthm}]
   \COMMENT{}

Let us show first that every non-empty cut $F$ with a small side $A$ contains a cocircuit of~$M$. If not, $F$~is independent in~$M^*$, so it avoids a base $B$ of~$M$. Adding an element $f\in F$ to $B$ creates a circuit of~$M$, a cycle or double ray that meets $F$ precisely in~$f$. Since $G[A]$ contains no ray, this is impossible.

Conversely, let us show that every cocircuit $D$ of $M$ contains a \xsided\ cut. Pick an edge $e\in D$. If its endvertices lie in the same component of $G-D$, then $G$ contains a cycle meeting $D$ in exactly~$e$, contradicting Lemma~\ref{circuitcocircuitcap}. So the endvertices of $e$ lie in distinct components of $G-D$. If both these contain a ray, then these rays can be chosen so as to combine with $e$ to a double ray meeting $D$ precisely in~$e$, again contradicting Lemma~\ref{circuitcocircuitcap}. Hence one of these components contains no ray. Its vertex set $A$ is the small side of a cut $F$ with $e\in F\sub D$.%
   \COMMENT{}
   To show that $F$ is a skew cut, we still have to show that it contains no non-empty cut $F'$ with a small side properly. But any such $F'$ contains a cocircuit~$D'$, as shown earlier, giving $D'\sub F'\subne F\sub D$. This contradicts the fact that no cocircuit contains another cocircuit properly.

We have shown that every skew cut contains a cocircuit, and vice versa. Since skew cuts, as cocircuits, cannot contain each other properly, these inclusions cannot be proper. So the cocircuits of $M$ are the skew cuts of~$G$.
\end{proof}

\section{Thin-sum matroids and representability}\label{sec:thinsums}

In finite matroid theory, representable matroids are an important generalization of graphic matroids. As matroids defined by linear independence are finitary, the dual of an infinite representable matroid will not, except in trivial cases, be representable. Representability, as usually defined, is thus another concept that seems too narrow for infinite matroids.

In this section we present a notion of vector independence, different from linear independence, that can give rise to non-finitary matroids. Examples include the algebraic cycle matroids of graphs and of higher-dimensional complexes~\cite{InfMatroidAxioms}.

Let $F$ be a field, and let $A$ be some set. We say that a set $X$ of 
functions $x: A\to F$ is \emph{thin\/} if for every $a\in A$ there are only finitely many $x\in X$ 
with $x(a)\neq 0$. Given such a thin set of functions, their pointwise sum $\sum_{x\in X}x$ 
is another $A\to F$ function.%
   \COMMENT{}
   We say that a set of $A\to F$ functions, not necessarily thin,%
   \footnote{Requiring independent sets to be thin leads to a different notion of representability that may have its own applications. It is easily seen that this notion does not satisfy (IM) for arbitrary sets $E$ of $A\to F$ functions, but there may be interesting examples where it does.}%
   \COMMENT{}
   is \emph{thinly independent}
if for every thin subset $X$ and every corresponding%
   \COMMENT{}
   family $(\alpha_x)_{x\in X}$ of coefficients $\alpha_x\in F$ 
we have $\sum_{x\in X}\alpha_x x={\bold 0}\in F^A$ only when $\alpha_x=0$ for all $x\in X$.%
   \COMMENT{}

Unlike linear independence, thin independence does not always define a matroid.%
   \footnote{View the elements of $E=\mathbb F_2^\mathbb N$ as subsets of $\mathbb N$,
and define sets $I:=\{\{1,n\}:n\in\mathbb N\}$ and $I':=\{\{n\}:n\in\mathbb N\}$.
Both $I$ and $I'$ are thinly independent. Moreover, $I'$ is maximally thinly independent
but $I$ is not: $I+\mathbb N$, for instance, is still thinly independent. 
Yet, the only $x\in I'$ for which $I+x$ is thinly independent is $x=\{1\}$,
which is already contained in $I$. Thus, (I3) is violated.}
   The following theorem gives a sufficient condition for when it does:%
   \COMMENT{}

\begin{theorem}\label{thinsumsmatroid}
If a set $E$ of $A\to F$ functions is thin, then its thinly independent subsets form
the independent sets of a matroid on~$E$.%
   \COMMENT{} 
\end{theorem}

Let us call such a matroid as in Theorem~\ref{thinsumsmatroid} the {\it thin-sums matroid\/} of the functions in~$E$. In the remainder of this section we prove Theorem~\ref{thinsumsmatroid}, and then briefly discuss what it means for graphs.

\bigskip
Let $F$ be endowed with the discrete topology,%
   \COMMENT{}
   and the set $F^A$ of all $A\to F$ functions with the product topology. Thus for each $x\in F^A$, the sets
 $$\{\,y\in F^A : y(a)=x(a) \text{ for all }a\in A'\}$$
where $A'$ ranges over the finite subsets of $A$ forms a basis of the open neighbourhoods of~$x$.%
   \COMMENT{}

   Given a set $X$ of functions $A\to F$, we write $\langle X\rangle$ for the set
of all functions $A\to F$  that are of the form $\sum_{x\in X'} \alpha_x x$, 
where $X'\subseteq X$ and $X'$ is thin.%
   \COMMENT{}
Similarly, we write $\overline X$ for the closure in~$F^A$ of the set $X$.

In contrast to~\cite{basis} where these concepts were introduced,
$X$ will here always be a subset of a thin set~$E$. While the sets $\langle X\rangle$ and $\overline X$ may contain elements outside our ground set $E$, we note that $X\mapsto \langle X\rangle\cap E$ is the closure operator associated with the set $\mathcal I$ of 
thinly independent subsets of $E$, as defined in Section~\ref{sec:axioms}.

We need two lemmas from~{\rm\cite{basis}}, which together imply that $\langle\langle X\rangle\rangle = \langle X\rangle$ for all our (thin) sets $X\sub E$:

\begin{lemma}{\rm\cite[Lemma~5]{basis}}\label{lem:closed}
Every thin set $X\sub F^A$ satisfies $\overline{\left<X\right>}=\left<X\right>$.
\end{lemma}

\begin{lemma}{\rm\cite[Lemma~6]{basis}}\label{spanspan}
Every set $X\sub F^A$ satisfies $\langle\langle X\rangle\rangle\subseteq \overline{\langle X\rangle}$.
\end{lemma}

\begin{proof}[Proof of Theorem~\ref{thinsumsmatroid}]
Let $\mathcal I$ be the set of thinly independent subsets of $E$.
   Clearly, $\I$ satisfies (I1) and (I2).%
   \footnote{We use the independence axioms in our proof. Alternatively, one could check that 
$\langle\,\cdot\,\rangle\cap E$ is indeed the closure operator associated with~$\I$ and then use the closure axioms: (CL1), (CL2) and (CL4) are straightforward, (CL3) follows from our two lemmas, and (CLM) is proved like (IM) in the text.}

Our first claim is that, for all sets $J\subseteq X\subseteq E$ with $J\in\I$,
\begin{equation}
\label{clm:1}
\emtext{if $X\subseteq\langle J\rangle$ then $J$ is a 
maximal element of $\{I\in\I:I\subseteq X\}$}.
\end{equation}
Consider an $x\in X\sm J$. Then there are coefficients
$a_j\in F$, $j\in J$, such that
$\sum_{j\in J} a_jj = x$. Thus, $J+x$ is not thinly independent for any $x\in X\sm J$, which 
implies the claim.

To prove that $\I$ satisfies~(I3), let $I\in\I\sm\Imax$ and
$I'\in\Imax$ be given. We have to find an $x\in I'\sm I$ such that $I+x$ is still thinly independent. Clearly, any $x$ in $I'\sm\langle I\rangle$ will do, so it suffices to show that $I'\not\subseteq\langle I\rangle$. If $I'\subseteq\langle I\rangle$, then $\langle I'\rangle \subseteq \langle\langle I\rangle\rangle=\langle I\rangle$ by Lemmas~\ref{lem:closed} and~\ref{spanspan}.
As $I'\in\Imax$ implies $E\subseteq\langle I'\rangle$,%
   \COMMENT{}
   this yields $E\subseteq\langle I'\rangle\subseteq\langle I\rangle$, which contradicts~\eqref{clm:1}. This completes the proof of~(I3). 

Next, we prove that $\I$ satisfies~(IM). This will follow directly from~\eqref{clm:1} and the following claim:
\begin{equation}
\label{clm:2}
\begin{minipage}[c]{0.8\textwidth}\em
For all sets $I\sub X\subseteq E$ with $I\in\I$, there is a $B\in\I$ 
with $I\subseteq B\subseteq X$ and $X\subseteq\langle B\rangle$.
\end{minipage}\ignorespacesafterend 
\end{equation} 

In the remainder of this proof we thus prove~\eqref{clm:2}. Let $x_1,x_2,\ldots$ be a (possibly transfinite) enumeration of $X\sm I$.
Inductively, we define nested sets $B_\lambda\subseteq X$
as follows. Start with $B_0:=I$.
If, in step~$\lambda$, there are families $(\alpha_\mu)_{\mu>\lambda}$ and $(\beta_i)_{i\in I}$
of coefficients in $F$ such that
 $$x_\lambda = \sum_{\mu>\lambda}\alpha_\mu x_\mu+\sum_{i\in I}\beta_i i$$
(these sums are well-defined, since $E$ is thin by assumption), we set 
$B_{\lambda}:=\bigcup_{\mu<\lambda}B_\mu$.
Otherwise we put 
$B_{\lambda}:=\{x_\lambda\}\cup \bigcup_{\mu<\lambda}B_\mu$.
Finally, we let 
$B:=\bigcup_{\lambda}B_\lambda$,
which we claim satisfies~\eqref{clm:2}.

Let us first check that $B\in\I$. If not then
there are coefficients $\alpha_b\in F$ for all $b\in B$, not all of them zero,
such that $\sum_{b\in B} \alpha_bb=\mathbf 0$.%
   \COMMENT{}
   Since $I$ is thinly independent, 
there must be some $b\in B\sm I$ with $\alpha_b\neq 0$.
Pick such a function $b=x_\lambda$ with smallest index~$\lambda$.
Then 
 $$-x_\lambda\ = \sum_{x\in \{x_\mu : \mu>\lambda\}\cap
B}\alpha^{-1}_{x_\lambda}\alpha_xx+\sum_{i\in I}\alpha^{-1}_{x_\lambda}\alpha_ii\,,$$
which contradicts the fact that $x_\lambda\in B$.

Next, we prove  $X\subseteq\overline{\left<B\right>}$.
By Lemma~\ref{lem:closed} this will imply 
$X\subseteq\overline{\left<B\right>}=
\left<B\right>$,%
   \COMMENT{}
   which then completes the proof of~\eqref{clm:2}.

To prove $X\subseteq\overline{\left<B\right>}$, consider a function $z\in X$. We need to find, for every finite subset $A'$ of~$A$, a function $z'\in\langle B\rangle$ that agrees
with $z$ on $A'$. 
Denote by $L$ the set of $x\in X$ for which there is an $a\in A'$ with $x(a)\neq 0$.
Observe that $L$ is a finite set, since $A'$ is finite and $X\subseteq E$ is thin. 
In particular, we may write $L\sm B=\{x_{\lambda_1},\ldots,x_{\lambda_k}\}$, 
where $\lambda_1<\ldots<\lambda_k$. Let $\ell$ be the largest number in $\{1,\ldots,k+1\}$ for which there exists a $y\in \langle \{x_{\lambda_\ell},\ldots,x_{\lambda_k}\}\cup(B\cap L)\rangle$
with $y(a)=z(a)$ for all $a\in A'$.
Note that there is always such a $y$ for $\ell=1$, as we may either pick 
$y=z$ if $z\in L$, or  $y=\mathbf 0$ otherwise. Note, furthermore, 
that we have found the desired $z'$ if $\ell=k+1$, as then $z':=y\in\langle B\rangle$. 

Suppose that $\ell\leq k$.
Since $x_{\lambda_\ell}\notin B$ there are coefficients 
$(\alpha_\mu)_{\mu>{\lambda_\ell}}$ and $(\beta_i)_{i\in I}$
such that $x_{\lambda_\ell}=\sum_{\mu>\lambda_\ell}\alpha_\mu x_\mu+\sum_{i\in I}\beta_ii$.
Then
  $$x_{\lambda_\ell}=\sum_{\mu>\lambda_\ell,\,
x_\mu\notin B}\alpha_\mu x_\mu+\sum_{b\in B}\beta'_bb$$
  with suitable coefficients $\beta'_b$. Restricting this to~$L$, set
  $$r:=\sum_{p=\ell+1}^k\alpha_{\lambda_p}x_{\lambda_p}
+\sum_{b\in B\cap L}\beta'_bb\,.$$
Then 
\[
r(a)=x_{\lambda_\ell}(a) \emtext{ for all }a\in A', \emtext{ and }
r\in \langle \{x_{\lambda_{\ell+1}},\ldots,x_{\lambda_k}\}\cup(B\cap L)\rangle.
\]

Next, by choice of $\ell$ there 
exists $y=\sum_{q=\ell}^k\gamma_{\lambda_q} x_{\lambda_q}+\sum_{b\in B\cap L}\delta_bb$
such that $y(a)=z(a)$ for all $a\in A'$. 
Replacing $x_{\lambda_\ell}$ in this sum with~$r$, we obtain
 $$y^*:= \gamma_{\lambda_\ell}r +
\sum_{q=\ell+1}^k\gamma_{\lambda_q} x_{\lambda_q}
+\sum_{b\in B\cap L}\delta_bb\ \in\ \langle \{x_{\lambda_{\ell+1}},\ldots,x_{\lambda_k}\}\cup(B\cap L)\rangle.$$
As $y^*$ agrees with $z$ on~$A'$, this contradicts the maximal choice of $\ell$.
\end{proof}

The algebraic cycle matroid of a graph $G=(V,E)$ can be represented as a thin-sums matroid over~$\mathbb F_2$, for any $G$ for which it is defined (cf.\ Theorem~\ref{HiggsBean}). Indeed, as in finite graphs we represent an edge $e=uv$ by the map $V\to\mathbb F_2$ assigning 1 to both $u$ and~$v$, and 0 to every other vertex. Then a set $F\sub E$ of edges becomes a set of $V\to\mathbb F_2$ functions, not necessarily thin, which is thinly independent if and only if $F$ contains no elementary algebraic cycle.%
   \COMMENT{}
   This example can be generalized to higher dimensions; see~\cite{InfMatroidAxioms} for algebraic cycle matroids of simplicial complexes.

We do not know whether the other non-finitary matroids we discussed in this paper are representable as thin-sum matroids, but suspect not. For finitary matroids, one would hope that `thin-sum' representability coincides with traditional representability, but we have no proof of this:

\begin{problem}
Is a finitary matroid representable as a thin-sums matroid if and only if it is representable in the usual sense?%
   \COMMENT{}
\end{problem}

We can show that the finite-cycle matroid $\mfc(G)$ of a graph $G$ is a thin-sums matroid if $G$ has finite chromatic number, but we do not know this for arbitrary~$G$.%
   \COMMENT{}

\section*{Acknowledgement}
We would like to thank Robin Christian for indicating the proof of implication (iii)$\to$(i) of Theorem~\ref{graphicduals}, which we had previously stated as a conjecture.

\bibliographystyle{amsplain}
\bibliography{collective} %included explicitly below due to hand-edited bits such as "in this volume"

\small
\vskip2mm plus 1fill\noindent
Version 12 March 2011.\\ Footnotes 5 and 6 were added later. Figure~1 was corrected in May 2012.

\end{document}